\documentclass[12pt]{article}
\usepackage{}
\usepackage{amsmath}

\usepackage{epsf}

\usepackage{indentfirst}
\usepackage{latexsym}
\usepackage{amssymb}
\usepackage[dvips]{graphicx}
\usepackage{flafter}
\usepackage{caption2}
\usepackage{textcomp}
\usepackage{supertabular}
\usepackage{longtable, booktabs}

\usepackage{multirow}

\usepackage{hhline}
\usepackage{bigstrut,bigdelim}

\usepackage{rotating}
\usepackage{multicol}
\usepackage{multirow}
\usepackage{makeidx} 
\usepackage{epsfig}
\usepackage{fancyhdr}       

\newif\ifger

\gerfalse

\newtheorem{theorem}{Theorem}[section]

\newtheorem{lemma}[theorem]{Lemma}

\newtheorem{corollary}[theorem]{Corollary}

\newtheorem{algorithm}[theorem]{Algorithm}

\def\C{{\mathfrak C}}
\def\B{\mathcal{B}}
\def\D{\mathcal{D}}\def\G1{G^{\C}}

\newcommand{\cB}{{\cal B}}

\newcommand{\cD}{{\cal D}}

\newcommand{\cP}{{\cal P}}

\def\la{{\langle}}

\def\lam{\lambda}

\def\Soc{{\rm Soc}}
\def\Aut{{\rm Aut}}

\newcommand{\Out}{{\rm Out}}

\newcommand{\PSL}{{\rm PSL}}
\newcommand{\PGL}{{\rm PGL}}
\newcommand{\PSSL}{{\rm P\Sigma L}}
\newcommand{\PGGL}{{\rm P\Gamma L}}

\def\la{\lambda}

\def\PSL{{\rm PSL}}
\def\PGL{{\rm PGL}}

\def\Soc{{\rm Soc}}

\begin{document}

\title{Classification of flag-transitive primitive symmetric $(v,k,\la)$ designs with  $\PSL(2,q)$ as socle}

\author{Shenglin Zhou$^{a}$\footnote{Supported by the National Natural Science Foundation
of China (No.11471123) and the Natural Science Foundation of
Guangdong Province (No.S2013010011928). E-mail address:
slzhou@scut.edu.cn. The results of this paper is part of Ph.D. thesis of the second author.}
, Delu Tian$^{a, b}$\\
{\small\it  a. School of Mathematics, South China University of Technology,}\\
 {\small\it Guangzhou, Guangdong 510640, P. R. China}\\
 {\small \it b. Department of
Mathematics, Guangdong University of
Education,}\\
{\small\it Guangzhou, 510310,  P. R. China}\\
\date{}
 }
\maketitle
\date{}

\begin{abstract}
Let $\cal D$ be a nontrivial symmetric  $(v,k,\la)$ design, and
$G$ be a subgroup of  the full automorphism group of $\cal D$.
In this paper we prove that if  $G$ acts flag-transitively, point-primitively on $\D$
and $\Soc(G)= \PSL(2,q)$, then $\cal D$ has parameters $(7, 3, 1)$,
$(7, 4, 2)$, $(11, 5, 2)$, $(11, 6, 3)$ or $(15, 8, 4)$.


\medskip
\noindent{\bf MR(2000) Subject Classification} 05B05, 05B25, 20B25

\end{abstract}

\maketitle

\section{Introduction}

 A  2-$(v,k,\la)$ {\it design} $\D$ is a set $P$ of $v$ points together with a set $\B$ of $b$ blocks,  such that
every block contains $k$ points and every pair of points is in
exactly $\la$ blocks. The design  $\D$  is  {\it symmetric} if
$b=v$, and is {\it non-trivial} if $2< k< v-1$. In this paper we
only study non-trivial symmetric 2-$(v, k, \lambda)$ designs, and
for brevity we call such a design a symmetric $(v, k, \lambda)$
design. A {\it flag} in a design is an incident point-block pair.
The {\it complement} of $\cal D$, denoted by $\cal D'$,  is a
symmetric $(v, v-k, v-2k+\lam)$
 design whose set of points is the same as the set of points of $\cal D$, and whose blocks
are the complements of the blocks of $\cal D$.
  The {\it automorphism group $\Aut(\cal D)$} of ${\cal D}$ consists of
all permutations of $P$ which leave $\cal B$ invariant. For $G
\leq \Aut(\cal D)$, the design $\cal D$ is called {\it
point-primitive} if $G$ is primitive on $P$, and {\it
flag-transitive} if $G$ is transitive on the set of flags. The
{\it socle} of a group $G$, denoted by $\Soc(G)$, is the subgroup
generated by its minimal normal subgroups.

The classification program for symmetric $(v,k,\la)$ designs has
been studied by several researchers. In 1985, Kantor
\cite{Kantor2-Tr} classified all symmetric $(v,k,\la)$ designs
admitting  2-transitive automorphism groups. In \cite{Dempwolff
Rank3}, Dempwolff  determined all symmetric $(v,k,\la)$ designs
which admit an automorphism group $G$
 such that $G$ has a nonabelian socle and is a primitive rank three group on points (and blocks).
 In \cite{TZ Sporadic}, we classified flag-transitive point-primitive symmetric  $(v,k,\lambda)$ designs
  admitting an automorphism group $G$ such that $\Soc(G)$ is a sporadic simple group.
 This paper is devoted to the complete classification of flag-transitive point-primitive symmetric $(v,k,\la)$ designs
 which admit an automorphism group $G$ with $\Soc(G)= \PSL(2,q)$,
 and extend the result of symmetric designs with $\lambda=4$ in \cite{ZT Classical4} to the
 general case.

\begin{theorem}\label{MainThm}
Let $\D=(P,\B)$ be a  symmetric $(v,k,\lam)$ design which admits a flag-transitive, point-primitive
 automorphism group $G$, and $x$ be a point of $P$. If $G$ is an almost simple group and
 $X=\Soc(G)=\PSL(2,q)$, where $q = p^f$ and $p$ is a prime,
 then $\D$ is one of the following:
  \begin{enumerate}
  \item[(i)] \,
        a $(7, 3, 1)$ design with  $X=\PSL(2,7)$ and $X_x=S_4$;
  \item[(ii)] \,
    a $(7, 4, 2)$  design with  $X=\PSL(2,7)$ and $X_x=S_4$;
  \item[(iii)] \,
    a $(11, 5, 2)$ design with  $X=\PSL(2,11)$ and $X_x=A_5$;
  \item[(iv)] \,
     a $(11, 6, 3)$ design with  $X=\PSL(2,11)$ and $X_x=A_5$;
  \item[(v)] \,
    a $(15, 8, 4)$ design with    $X=\PSL(2,9)$ and $X_x=\PGL(2,3)$.
  \end{enumerate}
\end{theorem}

\begin{corollary}\label{coro}
For $\la\geq 5$, there is no symmetric $(v,k,\lam)$ design
admitting a flag-transitive, point-primitive almost simple
 automorphism group with socle $\PSL(2,q)$.
\end{corollary}

\section{Preliminaries}

In this section we  state some  preliminary results which will be
needed later in this paper. From \cite{Reg ClassicalBi} and
\cite{Zhou ClassicalTri} we get the following:

 \begin{lemma} \label{divide} \quad
 Let $\cD$ be a symmetric $(v,k,\lam)$ design. Then the following hold:
  \begin{enumerate}
  \item[(i)] \, $k(k-1)=\lam (v-1)$, and in particular $k^2>v$;
  \item[(ii)] \,  $k\mid \lam d_i$, where $d_i$ is any  non-trivial subdegree of $G$;
  \item[(iii)] \,   $k \mid  |G_x|$ and $|G_x|^3 > |G|$, where $G_x$ is the stabilizer in $G$ of a point $x\in P$.
  \end{enumerate}
\end{lemma}

\begin{lemma} \label{flag} {\rm (\cite{FiniteGeometries})}\quad 
 Let $\cD=(P, \cB)$ be a symmetric $(v,k,\lam)$  design and $G\leq \Aut(\cal D)$.
 Suppose that $\gcd(k, \lam)=1$ and $G$ is 2-transitive on  $P$.  Then $G$ is flag-transitive.
\end{lemma}

\begin{lemma} \label{rankPB} {\rm (\cite{SymDesApproach})}\quad
  Let $\cD$ be a symmetric $(v,k,\lam)$  design and $G\leq \Aut(\cal D)$. Then
  \begin{enumerate}
  \item[(i)] \,
  $G$ has as many orbits on points as on blocks;
  \item[(ii)] \,
 if $G$ is a transitive automorphism group, then $G$ has the same rank whether considered as a permutation group
    on points or on blocks.   
   \end{enumerate}
\end{lemma}

\begin{lemma}  \label{ComplementFlagTr} \quad
Let $\cD=(P, \cB)$ be a symmetric $(v,k,\lam)$  design admitting a
flag-transitive automorphism group $G$. Then  $G$ is 2-transitive
on $P$ if and only if $G$ is flag-transitive on
 $\cal D'$, the complement  design of $\cal D$.
\end{lemma}
 \textbf{Proof.}
 Suppose that $G$ is 2-transitive on $P$.
Then for any $x\in P$ and $B\in \cB$, Lemma \ref{rankPB} shows that both $G_x$ and $G_B$ acting on points or on blocks  have two orbits.
The flag-transitivity implies that $G_B$ acts transitively  on the
points of $B$. Thus $G_B$ has an orbit $\Gamma_1$ of length $k$ on
$P$ and the other orbit $\Gamma_2$ of length $v-k$ and
$\Gamma_2=P-\Gamma_1$. Therefore, $\Gamma_2 = B'$ is one of the
blocks of $\cD'$ and  $G_B=G_{B'}$. So $G_{B'}$ is transitive on
the points of $B'$. Moreover, by Lemma \ref{rankPB},
 $G$ is block-transitive on $\cD'$ since
 $G$ is 2-transitive on $P$.
Hence  $G$ is flag-transitive on $\cD'$.

Conversely, if $\cal D'$  is flag-transitive, then $G_{B'}$ is
transitive on the points of $B'$ for every block $B'$ of $\cal
D'$. Let $B=P-B'$. Then $B$ is one of the blocks of $\cD$ and
$G_B=G_{B'}$. Since $\cal D$  is flag-transitive,  $G_{B}$ is
transitive on $B$. Thus $G_B$ has two orbits acting on points,
which implies that the point-stabilizer $G_x$ has two orbits
acting on $P$ by Lemma \ref{rankPB}. Hence $G$ is 2-transitive on
$P$.                     $\hfill\square$

\begin{lemma}  \label{orbits} \quad
Let $G$ be a group that acts transitively on $P$, and let $X\unlhd G$.
Then  each orbit  of $G_x$,  the stabilizer of $G$ for some $x\in P$, acting on $P$ is the union of some
orbits of $X_x$ which have the same  cardinality.
\end{lemma}
 \textbf{Proof.}
 Suppose that
 the orbits of the action of the stabilizer $X_x$ on $P$ are $\Delta _1(=\{ x\})$, $\Delta _2, \cdots$, $\Delta _m$.
 For any  orbit $\Delta _i$ ($i=1,2,\cdots,m$) and any $g\in G_x$,
 we have $(\Delta _i^g)^{X_x}=(\Delta _i^g)^{g^{-1}{X_x}g}=\Delta _i^{X_xg}=\Delta _i^g$,
 where the first equality holds since $X_x\unlhd G_x$.
 Hence $\Delta_i^g$ is also an orbit of $X_x$.
 It follows that $\Delta_i^{G_x}=\{\Delta_i^g\ |\ g\in {G_x}\}$, a union of orbits of $X_x$  which have the same  cardinality,
 gives an orbit of $G_x$.
    $\hfill\square$

\begin{lemma}  \label{divide1} \quad
Let $\D=(P,\B)$ be a symmetric $(v,k,\la)$ design admitting a flag-transitive, point-primitive automorphism group $G$ with socle $X$.
 If the non-trivial subdegree $t$ of $X$ appears with multiplicity $s$, then $k\mid \la st$.
\end{lemma}
 \textbf{Proof.}
Suppose that $\Gamma_1$, $\Gamma_2, \cdots$, $\Gamma_s$ are all orbits of $X_x$ with cardinality $t$, where $x\in P$.
By Lemma \ref{orbits}, the group $G_x$ acts on $\Gamma=\bigcup \limits^{s}_{i=1}\Gamma_i$,  and  the cardinalities of orbits of $G_x$ are
$$(a_1t)^{s_1},\ (a_2t)^{s_2},\ \cdots, \  (a_rt)^{s_r},$$
where $a^b$ means that $a$ appears with multiplicity $b$, and $r$, $a_i$, $s_i$ ($1\leq i\leq r$) are all positive integers such that
$\sum\limits^{r}_{i=1}a_{i}{s_i}=s$, and $a_i\ne a_j$ if and only if $i\ne j$ for $1\leq i, j\leq r$.
 Let $c=\gcd(a_1,a_2,\cdots,a_r)$, then $c\mid s$.
Lemma \ref{divide} ($ii$) shows that $k\mid  \lambda (a_it)$,
$i=1,2,\cdots,r$.
 So $k\mid \lambda ct$, and hence $k\mid \lambda st$.
    $\hfill\square$\\

The subgroups of $\PSL(2,q)$  are well-known and given by Huppert
\cite{Huppert-Dickson}.

\begin{lemma}  \label{AllSubgroups}{\rm (\cite{Huppert-Dickson})}\quad
 The subgroups of the group $\PSL(2,q)$ $(q=p^f)$ are as follows.
  \begin{enumerate}
\item[(i)] \,
 An elementary abelian group  $C_{p}^{\ell}$, where $\ell \leq f$.
\item[(ii)] \,
 A cyclic group $C_z$,  where $z\mid \frac{p^f\pm 1}{d}$ and $d=\gcd(2,q-1)$.
\item[(iii)] \,
 A dihedral group $D_{2z}$, where $z$ is the same as in (ii). 
\item[(iv)] \,
 The alternating group $A_4$ when $p>2$ or $p=2$ and $2\mid f$.
\item[(v)] \,
 The symmetric group $S_4$ when $p^{2f}\equiv 1\pmod {16}$.
\item[(vi)] \,
 The alternating group $A_5$ when $p=5$ or $p^{2f}\equiv 1\pmod {5}$.
\item[(vii)] \,
 $C_{p}^{\ell}:C_t$, where $t\mid \gcd(p^\ell-1, \frac{p^f-1}{d})$ and $d=\gcd(2,q-1)$.
\item[(viii)] \,
 $\PSL(2,p^\ell)$ when $\ell\mid f$ and $\PGL(2,p^\ell)$ when $2\ell\mid f$.
   \end{enumerate}
\end{lemma}

The following lemma is a combination of Theorem 1.1, 2.1 and 2.2 in \cite{MaximalSubgroupsOfPSL(2q)}.

\begin{lemma} \label{MaxSubgp} \quad
 Let $X = \PSL(2,q)\leq G \leq \PGGL(2,q)$ and let $M$ be a
maximal subgroup of $G$ which does not contain $X$. Then either
$M\cap X$ is maximal in $X$, or $G$ and $M$ are given in Table
\ref{tab:1}. The maximal subgroups of $X$ appear in Tables
\ref{tab:2} and \ref{tab:3}.
\begin{table}[h]
\begin{center}
\begin{tabular}{llc}
\hline
$G$&$M$ &$|G:M|$\\
\hline
$\PGL(2,7)$ & $N_G(D_6) = D_{12}$&28\\
$\PGL(2,7)$ &$N_G(D_8) = D_{16}$&21\\
$\PGL(2,9)$ & $N_G(D_{10}) = D_{20}$&36\\
$\PGL(2,9)$ &$N_G(D_8) = D_{16}$&45\\
 $M_{10}$ &$N_G(D_{10}) = C_5\rtimes C_4$&36\\
$M_{10}$ & $N_G(D_8) = C_8 \rtimes C_2$& 45\\
$\PGGL(2,9)$ & $N_G(D_{10}) = C_{10}\rtimes C_4$&36\\
 $\PGGL(2,9)$  & $N_G(D_8)= C_{8}.\Aut(C_8)$&45\\
$\PGL(2,11)$ &$N_G(D_{10}) = D_{20}$&66\\
$\PGL(2,q), q = p\equiv \pm 11, 19({\rm mod} \ {40})$ & $N_G(A_4)=S_4$& $\frac{q(q^2-1)}{24}$\\
\hline
\end{tabular}
\caption{$G$ and $M$ of Lemma \ref{MaxSubgp}} \label{tab:1}
\end{center}
\end{table}
\begin{table}[h]
\begin{center}
\begin{tabular}{llcc}
\hline
   Structure & Conditions  & Order & Index\\
\hline
   $C_p^f:C_{(q-1)/2}$ &             &  $\frac{q(q-1)}{2}$  & $q+1$\\
  $D_{q-1}$           & $q\geq 13$  &  $q+1$               & $\frac{q(q-1)}{2}$ \\
                $D_{q+1}$           & $q\ne 7,9$  &  $q-1$               & $\frac{q(q+1)}{2}$ \\
                $\PGL(2,q_0)$        & $q=q_0^2$   &  $q_0(q_0^2-1)$      & $\frac{q_0(q_0^2+1)}{2}$ \\
                $\PSL(2,q_0)$        & $q=q_0^r$, $r$ odd  prime   &  $\frac{q_0(q_0^2-1)}{2}$      & $\frac{q_0^{r-1}(q_0^{2r}-1)}{q_0^2-1}$ \\
                $A_5$               & $q=p\equiv \pm 1({\rm mod} \ {5})$,   &  120      &  $\frac{q(q^2-1)}{120}$ \\
                                    & or $q=p^2\equiv -1({\rm mod} \ {5})$   &      &  \\
                $A_4$               & $q=p\equiv \pm 3({\rm mod} \ {8})$,   &  12      & $\frac{q(q^2-1)}{24}$ \\
                                    & and $q\not\equiv \pm 1({\rm mod} \ {10})$   &      &  \\
                $S_4$               & $q=p\equiv \pm 1({\rm mod} \ {8})$   &  24      &
                $\frac{q(q^2-1)}{24}$\\
\hline
\end{tabular}
 \caption{Maximal subgroups of $\PSL(2,q)$} with $q=p^f\geq 5$, $p$ odd prime \label{tab:2}
\end{center}
\end{table}
\begin{table}[h]
\begin{center}
\begin{tabular}{llcc}
\hline
    Structure & Conditions  & Order & Index\\
\hline
 $C_2^f: C_{q-1}$         &             &  ${q(q-1)}$          & $q+1$ \\
                $D_{2(q-1)}$        &             &  $2(q+1)$            & $\frac{q(q-1)}{2}$ \\
                $D_{2(q+1)}$        &             &  $2(q-1)$            & $\frac{q(q+1)}{2}$ \\
                $\PSL(2,q_0)$        & $q=q_0^r$, $r$ prime, $q_0\ne 2$   &  $q_0(q_0^2-1)$      &
                $\frac{q_0^{r-1}(q_0^{2r}-1)}{q_0^2-1}$\\
\hline
\end{tabular}
 \caption{Maximal subgroups of $\PSL(2,q)$} with $q=2^f\geq 4$ \label{tab:3}
\end{center}
\end{table}
\end{lemma}
 Now we state the following algorithm, which will be useful to search for symmetric designs which satisfy the condition ``$k\mid u$''.
 The output of the algorithm is the list {\sc Designs} of parameter sequences $(v,k,\lam)$ of potential symmetric designs.

 \begin{algorithm}  \label{algorithm} \quad {\rm
 ({\sc Designs})\\
 {\sc Input}:$~~~\quad$    $u$, $v$.\\
 {\sc Output}:$~~~$  The list {\sc Designs} := $S$.\\
 set $S$ := an empty list;\\
 for each $k$ dividing $u$ with $2< k< v-1$\\
  $~~~~~~~ \lam:=k*(k-1)/(v-1)$;\\
  $~~~~~~$ if $\lam$ be an integer\\
  $~~~~~~~~~~$ Add $(v,k,\lam)$ to the list $S$;\\
  return $S$.
 }
 \end{algorithm}

\section{Proof of Theorem 1.1}

Let $\D=(P,\B)$ be a symmetric $(v,k,\la)$ design admitting a flag-transitive, point-primitive automorphism group $G$ with
 $X\unlhd G\leq \Aut(X)$, where $X=\PSL(2,q)$ with $q=p^f$ and $p$ prime.
As a maximal subgroup of $G$, the point stabilizer $G_x$ does not contain  $X$ since $X$ is transitive on $P$. 
 Thus Lemma \ref{MaxSubgp} shows that  either $X\cap G_x$ is maximal in $X$, or
$G$ and $G_x$ are given in Table 1. We will prove Theorem \ref{MainThm} by the
following three subsections.

\subsection{Cases in Table \ref{tab:1}}

In these cases, we may view  the maximal subgroup $M$ as the point stabilizer  $G_x$.
 We get the 3-tuples $(|G|, u, v)$ in Table \ref{tab:1}
  where $v$ is  the index $|G:G_x|$ and $u=|G_x|$. For each case except the last one,  we can obtain all potential symmetric designs using
 Algorithm \ref{algorithm} implemented in {\sc GAP} \cite{GAP4}. 
 There exists only one potential  $(21, 16, 12)$ design with $G=\PGL(2,7)$ and $G_x=D_{16}$.
 The subdegrees of $\PGL(2,7)$ acting on the cosets of  $D_{16}$  are 1, 4, 8 and
 8. (Throughout this paper, we apply {\sc Magma} \cite{Magma} to calculate the subdegrees of $G$ and the number of the conjugacy class of
 subgroups.)
 Then by using the {\sc Magma}-command
 {\tt Subgroups(G:OrderEqual:=n)} where $n=|G|/v$, we obtain the fact that $G$ has only one conjugacy class of subgroups with index 21.
  Thus $G_x$ is conjugate to $G_B$ for any $x\in \cP, B\in \cB$ which forces that there exists a block $B_0$ such that $G_x=G_{B_0}$.
 The flag-transitivity of $G$ implies that $G_{B_0}$ is transitive on the block $B_0$.
 So $B_0$ should be an orbit of $G_x$,
 but there is no such orbit of size $k=16$, a contradiction.

 Now we consider the last case. Here $G=\PGL(2,q)$ with $q=p\equiv \pm 11,19 \pmod {40}$,  $G_x=S_4$, and $v=\frac{q(q^2-1)}{24}$.
 Since $|G_x|^3>|G|$, we have $24^3>q(q^2-1)$, and so $q=p=11$ or 19.
 If $q=11$ then $v=55$.
 There exist two potential symmetric designs with parameters  $(55, 27, 13)$ and $(55, 28, 14)$,
  but neither of them satisfies the condition that $k\mid  |G_x|$.
 If $q=19$ then $v=285$, and so $(k,\la)=(72, 18)$ or $(213,
 159)$. However, for every case $k> 24$ which contracts the fact that $k$ divides $|G_x|$.

\subsection{Odd characteristic}

In this subsection, we consider the cases that $G$ has odd characteristic $p$ and $X\cap G_x$ is maximal in $X$.
The structure of $X\cap G_x$ comes from Table \ref{tab:2}.

\subsubsection*{{\bf Case (1)}. $X\cap G_x= E_q:\frac{q-1}{2}$.}

Here $v=q+1$, so $k(k-1)=\la (v-1)=\la q=\la p^f$. If $p\mid  k$,
then from $\gcd(p,k-1)=1$ we have $p^f\mid  k$, that is, $v-1 \mid
k$ which is impossible. Then $p\nmid k$, and so $p^f\mid k-1$
implies $v-1\mid  k-1$, which contradicts $k<v-1$.

 \subsubsection*{{\bf Case (2)}. $X\cap G_x= D_{q-1}$ $(q\geq 13)$.}  

In this case,  $v=\frac{1}{2}q(q+1)$,   $|\Out(X)|=2f$,
$|G|=\frac{1}{2}eq(q^2-1)$ and $|G_x|=e(q-1)$, where $e$ is a
positive integer and  $e\mid  2f$.

From $k\mid  |G_x|$  we get $k\mid  e(q-1)$. So there exists a
positive integer $m$ such that $k=\frac{e(q-1)}{m}$. The equality
$k(k-1)=\lam (v-1)$ implies that
$\frac{e(q-1)}{m}\big(\frac{e(q-1)}{m}-1\big)=\frac{1}{2}\lam
(q+2)(q-1)$, and hence
$$(2e^2-m^2\lam)q=2m^2\lam+2e^2+2em>0.$$
This implies that $m<2e$.
From $ p^f=q=\frac{2m^2\lam+2e^2+2em}{2e^2-m^2\lam}=\frac{6e^2+2em}{2e^2-m^2\lam}-2$, we get
$$p^f<6e^2+2em<6e^2+(2e)^2=10e^2\leq 40f^2.$$
It follows that the  3-tuples $(q, p, f)$ are
$$
\begin{array}{llllllll}
 &(27, 3, 3), &(81, 3, 4), &(243, 3, 5), &(729, 3, 6), &(25, 5, 2), &(125, 5, 3),\\
 & (625, 5, 4), &(49, 7, 2), &(343, 7, 3), &(121, 11, 2), &(13, 13, 1), &(17, 17, 1),\\
  &(19, 19, 1), &(23, 23, 1),  &(29, 29, 1), &(31, 31, 1),& (37, 37,
  1).&&\\
\end{array}
$$

  We call each of these  3-tuples  a subcase.
Since $k\mid  e(q-1)$ and $e\mid 2f$, it follows that $k\mid u$,
where $u=2f(q-1)$. It is easy to compute the values of $u$ and $v$
for every subcase. However, for every subcase, there is no such
symmetric design satisfying the condition that $k\mid u$  by
Algorithm \ref{algorithm}
 calculated with {\sc GAP}.

\subsubsection*{{\bf Case (3)}. $X\cap G_x= D_{q+1}$ ($q\neq 7, 9$).}

Now $v=\frac{1}{2}q(q-1)$,   $|\Out(X)|=2f$,
$|G|=\frac{1}{2}eq(q^2-1)$, and $|G_x|=e(q+1)$, where $e$ is a
positive integer and  $e\mid  2f$.

Since $k\mid  |G_x|= e(q+1)$, it follows that there exists a
positive integer $m$ such that $k=\frac{e(q+1)}{m}$. Then
$\frac{e(q+1)}{m}\big(\frac{e(q+1)}{m}-1\big)=\frac{1}{2}\lam
(q-2)(q+1)$. So we have
$$(m^2\lam-2e^2)q=2e^2-2em+2m^2\lam=2(e-\frac{1}{2}m)^2+(2\lam-\frac{1}{2})m^2>0.$$
Thus $
p^f=q=\frac{2e^2-2em+2m^2\lam}{m^2\lam-2e^2}=\frac{6e^2-2em}{m^2\lam-2e^2}+2$
which gives
$$p^f<6e^2+2\leq 24f^2+2.$$
Combining this with $q\neq 7, 9$, we obtain all possible 3-tuples $(q, p, f)$:
$$
\begin{array}{lllllllll}
 &(27, 3, 3), &(81, 3, 4), &(243, 3, 5), &(729, 3, 6), &(5, 5, 1),&(25, 5, 2), &(125, 5, 3), \\
  &(49, 7, 2), &(11, 11, 1), &(13, 13,1),&(17, 17, 1), &(19, 19, 1), &(23, 23, 1).&\\
\end{array}
$$

Since $k\mid e(q+1)$ and $e\mid 2f$,  then $k\mid u=2f(q+1)$.
 The values of $v$ and $u$ can be calculated easily for each 3-tuple $(p,q,f)$.
 In fact, we get no such symmetric design  satisfying $k\mid u$  by Algorithm \ref{algorithm}
 calculated with {\sc GAP}.

\subsubsection*{{\bf Case (4)}. $X\cap G_x= \PGL(2, q^{\frac{1}{2}})=\PGL(2, q_0)$.}

Here $v=\frac{q_0(q_0^2+1)}{2}$, $|X_x|=|X\cap G_x|=q_0(q_0^2-1)$,
$|\Out(X)|=2f$, $|G|=\frac{1}{2}eq(q^2-1)$, and
$|G_x|=eq_0(q_0^2-1)$, where $e\mid  2f$ and $f$ is even.

The subdegrees of $\PSL(2,q)$ on the cosets of $\PGL(2,q_0)$ are
       $$1,\ \frac{{q_0}({q_0}-\varepsilon)}{2},\ \ {q_0^2-1},\ ({q_0}({q_0}-1))^{\frac{{q_0}-4-\varepsilon}{4}},
       \ ({q_0}({q_0}+1))^{\frac{{q_0}-2+\varepsilon}{4}},$$
where ${q_0}\equiv \varepsilon~({\rm mod} \ {4})$ with
$\varepsilon =\pm 1$ (see \cite{Subdegree}). Recall that here
$a^b$ means the subdegree $a$ appears with multiplicity $b$. We
consider two subcases in the following.

\medskip
{\bf Subcase (4.1): }  $\varepsilon =-1$.  Then there exists a
positive integer $s$ such that $q_0=4s-1$, and the subdegrees here
are:
$$1,\ \frac{{q_0}({q_0}+1)}{2},\ {q_0^2-1},\
({q_0}({q_0}-1))^{\frac{{q_0}-3}{4}},\
({q_0}({q_0}+1))^{\frac{{q_0}-3}{4}}.$$ By Lemma \ref{divide1},
we get
 $$k\mid   \lambda\gcd\big(\frac{{q_0}({q_0}+1)}{2}, {q_0^2-1}, \frac{{q_0}({q_0}-1)({q_0}-3)}{4},
      \frac{{q_0}({q_0}+1)({q_0}-3)}{4}\big).$$
Since $q_0=4s-1$, it follows that $\gcd(\frac{{q_0}({q_0}+1)}{2},{q_0^2-1})=\frac{{q_0}+1}{2}$
  and  $\gcd(\frac{{q_0}({q_0}-1)({q_0}-3)}{4},
  \frac{{q_0}({q_0}+1)({q_0}-3)}{4})$ $=\frac{{q_0}(q_0-3)}{2}$.
 Thus $ k\mid \lambda\gcd(\frac{q_0+1}{2},  \frac{{q_0}(q_0-3)}{2})
    =2\lambda$.
 Then from $k>\la$ we get $k=2\la$.
By Lemma \ref{divide} ($i$),
$\la=\frac{v+1}{4}=\frac{q_0^3+q_0+2}{8}$ and
$k=\frac{q_0^3+q_0+2}{4}$. Since  $k \mid  |G_x|$ and $e\mid 2f$,
$k\mid u=2fq_0(q_0^2-1)$. It follows that $q_0^3+q_0+2\mid
8fq_0(q_0^2-1)$. Note that $\gcd(q_0^3+q_0+2, q_0)=1$ and
$\gcd(q_0^2-q_0+2, q_0-1)=2$, we get $q_0^2-q_0+2\mid 16f$. So $\
q_0^2-q_0+2\leq 16f$,  i.e.
$(p^{\frac{1}{2}f})^2-p^{\frac{1}{2}f}+2\leq 16f$. It follows that
$(f,p)=(2, 3)$ or $(2, 5)$ because $f$ is even. If $(f,p)=(2,5)$,
then $q_0$ is equal to $p$ which contradicts $q_0=4s-1$. Suppose
that $(f,p)=(2,3)$. Then $\cal D$ has parameters (15, 8, 4) with
$X=\PSL(2,9)\cong A_6$, $X_x=\PGL(2,3)\cong S_4$. The existence of
this design has been discussed in \cite{ZT Classical4}.

\medskip
{\bf Subcase (4.2):}  $\varepsilon =1$.  Then $q_0=4s+1$ for some
positive integer $s$.  Let $q_0=p^a$. Then $f=2a$. The subdegrees
are:
      $$1,\ \frac{{q_0}({q_0}-1)}{2},\ {q_0^2-1},\ ({q_0}({q_0}-1))^{\frac{{q_0}-5}{4}},\ ({q_0}({q_0}+1))^{\frac{q_0-1}{4}}.$$
 Lemma \ref{divide1} shows that
  $$k \mid  \lambda\gcd\big( \frac{{q_0}({q_0}-1)}{2}, {q_0^2-1},
 \frac{{q_0}({q_0}-1)({q_0}-5)}{4},
\frac{{q_0}({q_0}+1)(q_0-1)}{4}\big).$$ Since
$\gcd\big(\frac{{q_0}({q_0}-1)}{2}$, $ {q_0^2-1}\big)$ $=
\frac{1}{2}(q_0-1)$, it follows that $k \mid
\frac{1}{2}\lambda(q_0-1)$. Combining this with $k(k-1)=\lambda
(v-1)
 =\frac{1}{2}\lambda (q_0-1)(q_0^2+q_0+2)$, we get  $$q_0^2+q_0+2\mid k-1,$$ which implies that  $k$ is odd.

The flag-transitivity of $G$ implies that $G_x$ acts transitively
on $P(x)$,
 the set of all  blocks which are incident with the point $x$.
Therefore $G_x$ has some subgroup $L$ with index $k$. Since
$X_x\unlhd G_x$, we have $L/(L\cap X_x)\cong LX_x/ X_x$. Let
$H=L\cap X_x$, and $|LX_x:X_x|=c$ for some integer $c$. Then
$c\mid e$ and $|H|=\frac{eq_0(q_0^2-1)}{ck}$, and hence
$$k=\frac{e_0q_0(q_0^2-1)}{|H|},$$
 where $e_0=\frac{e}{c}$. The fact
$e\mid 2f=4a$ yields $e_0\mid 4a$.

 Since $\PSL(2,q_0)$ is the
normal subgroup of $\PGL(2,q_0)$ with index 2, and $H\leq
X_x=\PGL(2,q_0)$, we get $|H: H\cap \PSL(2,q_0)|
=|\PSL(2,q_0)H:\PSL(2,q_0)|=$ 1 or 2. Lemma \ref{AllSubgroups}
gives  all the subgroups of $\PSL(2,q_0)$, and hence $|H|$ must be
one of the following:

  \begin{enumerate}
  \item[(i)] \,
 $p^\ell$ or $2p^\ell$, where $\ell\leq a$;
  \item[(ii)] \,
 $z$ or $2z$, where $z\mid \frac{q_0\pm 1}{2}$;
  \item[(iii)] \,
 $2z$ or $4z$, where $ z\mid \frac{q_0\pm 1}{2}$;
  \item[(iv)] \,
 12 or 24;
  \item[(v)] \,
 24 or 48  when $p^{2a}\equiv 1\pmod {16}$;
  \item[(vi)] \,
 60 or 120  when $p=5$ or $p^{2a}\equiv 1\pmod {5}$;
  \item[(vii)] \,
 $tp^\ell $ or  $2tp^\ell $, where  $t\mid \gcd(p^\ell-1, \frac{p^a-1}{2})$;
  \item[(viii)] \,
 $\frac{1}{2}p^\ell(p^{2\ell}-1)$ or $ p^\ell(p^{2\ell}-1)$ when $\ell\mid a$,
   and  $p^\ell(p^{2\ell}-1)$ or $ 2p^\ell(p^{2\ell}-1)$ when $2\ell\mid a$.
     \end{enumerate}

Recall that $q_0=4s+1$, and so $8\mid q_0^2-1$. Combing this with
the  fact that $k=\frac{e_0q_0(q_0^2-1)}{|H|}$ is odd, gives
$8\mid |H|$. It follows that $|H|\neq p^\ell, 2p^\ell, z$, 12 and
60, and we deal with the remaining possible values of $|H|$ in
turn.

If $|H|=2z$ where $z\mid \frac{q_0\pm 1}{2}$ as in (ii) or (iii),
then it is easily known from $k=\frac{e_0q_0(q_0^2-1)}{2z}$ that
$k$ is even, a contradiction.

If $|H|=4z$ as in (iii), and in addition $z\mid \frac{q_0+1}{2}$,
then $k$ is even, a contradiction. Next suppose that  $z\mid \frac{q_0-1}{2}$. 
Since $q_0^2+q_0+2\mid k-1$,  then $q_0^2+q_0+2$ divides
$$
z(k-1)
=\frac{e_0q_0(q_0^2-1)}{4}-z=\frac{e_0(q_0-1)}{4}(q_0^2+q_0+2)-\frac{e_0(q_0-1)}{2}-z,$$
which implies that $q_0^2+q_0+2\mid \frac{e_0(q_0-1)}{2}+z$.
Therefore $q_0^2+q_0+2\leq \frac{e_0(q_0-1)}{2}+z$. It follows
that $p^{2a}+p^a+2 \leq 2a(p^a-1)+\frac{p^a-1}{2}$ because
$e_0\leq 4a$ and $z\leq \frac{q_0-1}{2}$, and hence
$2p^{2a}+2a+5\leq (2a-1)p^a$ which is a contradiction.

If $|H|=24$ then $k=\frac{e_0q_0(q_0^2-1)}{24}$. The fact that
$q_0^2+q_0+2\mid k-1$ implies that $q_0^2+q_0+2$ divides
$$
6(k-1)=\frac{e_0q_0(q_0^2-1)}{4}-6
=\frac{e_0(q_0-1)}{4}(q_0^2+q_0+2)-\frac{e_0(q_0-1)}{2}-6.$$ Thus
$q_0^2+q_0+2\mid \frac{e_0(q_0-1)}{2}+6$, and so $q_0^2+q_0+2\leq
\frac{e_0(q_0-1)}{2}+6$. Since $e_0\mid4a$, we have
$p^{2a}+p^a+2\leq \frac{e_0(q_0-1)}{2}+6 \leq 2a(p^a-1)+6$, which
is impossible since $p\geq 3$.

Now we turn to (v).  Here $|H|=48$ and
$k=\frac{e_0q_0(q_0^2-1)}{48}$. By $q_0^2+q_0+2\mid k-1$ we know
that $q_0^2+q_0+2$ divides
$$12(k-1)=\frac{e_0(q_0-1)}{4}(q_0^2+q_0+2)-\frac{e_0(q_0-1)}{2}-12,$$
which implies that $q_0^2+q_0+2\leq \frac{e_0(q_0-1)}{2}+12$.
Since $e_0\mid4a$, we have  $p^{2a}+p^a+2\leq
\frac{e_0(q_0-1)}{2}+12 \leq 2a(p^a-1)+12$,
which implies $p=3$ and $a=1$. Thus $q_0=3$, contradicting
$q_0=4s+1$.

 If $|H|=120$ then $k=\frac{e_0q_0(q_0^2-1)}{120}$. Using the fact
$q_0^2+q_0+2\mid k-1$,  we have that  $q_0^2+q_0+2$ divides
$$
30(k-1)
=\frac{e_0(q_0-1)}{4}(q_0^2+q_0+2)-\frac{e_0(q_0-1)}{2}-30.$$ So
$p^{2a}+p^a+2\leq \frac{e_0(q_0-1)}{2}+30 \leq 2a(p^a-1)+30$. It
follows that $(p,a)=(3,1)$ or $(5,1)$. The fact $q_0=4s+1$ shows
that $q_0=p=5$ and $a=1$. Now $e_0\mid 4a$ forces $e_0=1$, 2 or 4.
Hence $v=65$, and $k=1$, 2 or 4, respectively. However, $k$ is too small to satisfy $k^2>v$, a contradiction.

For (vii), if $|H|=tp^\ell$ or $2tp^\ell$ , then
$k=\frac{e_0q_0(q_0^2-1)}{itp^\ell}$ where $i=1$ or 2, and hence
$k$ is even because $t\mid \frac{p^a-1}{2}$, a contradiction.

For (viii), suppose first that $\ell\mid a$ and
$|H|=p^\ell(p^{2\ell}-1)$ or $\frac{1}{2}p^\ell(p^{2\ell}-1)$.
Then $k=\frac{ie_0p^a(p^{2a}-1)}{p^\ell (p^{2\ell}-1)}$ where
$i=1$ or 2. If $\ell=a$, then $k=ie_0$. From $v< k^2$ and $e_0\mid
4a$, we see $\frac{p^a(p^{2a}+1)}{2}<(ie_0)^2\leq 16i^2a^2$. It
follows that $(p,a)=(3,1)$, and so $q_0=3$, contradicting
$q_0=4s+1$. Thus $\ell<a$, and so $a\geq 2$. It is easy to see
that $p^\ell-1\mid p^a-1$ because $\ell \mid a$. Since
$q_0^2+q_0+2\mid k-1$, we obtain that $q_0^2+q_0+2$ divides
\begin{align*}
p^\ell(p^\ell+1)(k-1) & =\frac{ie_0p^a(p^{2a}-1)}{p^{\ell}-1}-p^\ell (p^\ell+1)\\
                      & =\frac{ie_0(p^a-1)}{p^{\ell}-1}(p^{2a}+p^a+2)-\frac{2ie_0(p^a-1)}{p^{\ell}-1}-p^\ell(p^\ell+1).
\end{align*}
Thus $p^{2a}+p^a+2\mid
\frac{2ie_0(p^a-1)}{p^{\ell}-1}+p^\ell(p^\ell+1)$. Since $\ell
\mid a$ and $\ell<a$, we have $2\ell\leq a$, and so $
p^{2\ell}\leq p^a$. Then $p^{2a}+p^a+2< 8ia(p^a-1)+2p^a$.
Combining this with $q_0=p^a=4s+1$ and $a\geq 2$, gives $(p,a)=(3,
2)$ when $i=1$, and $(p,a)=(3, 2)$ or (5, 2) when $i=2$. It
follows that $\ell=1$ and $e_0$=1, 2, 4 or 8. For all these
parameters $e_0, p, a$ and $\ell$, we can get all possible values
of $v$ and $k$. It is not hard to check that for all these pairs
$(v,k)$, there are no integer values of $\la$ satisfying equation
$k(k-1)=\la(v-1)$, a contradiction.

Now suppose that $2\ell\mid a$ and $|H|=p^\ell(p^{2\ell}-1)$. Then
$a\geq 2$ and $k=\frac{e_0p^a(p^a+1)(p^a-1)}{p^\ell
(p^{2\ell}-1)}$ is even since $p^{2\ell}-1\mid p^a-1$, a
contradiction.
 Finally suppose that $2\ell\mid a$ and
$|H|=2p^\ell(p^{2\ell}-1)$ so that
$k=\frac{e_0p^a(p^{2a}-1)}{2p^\ell (p^{2\ell}-1)}$ and $a\geq 2$.
Then by $q_0^2+q_0+2\mid k-1$, we get that $q_0^2+q_0+2$ divides
\begin{align*}
p^\ell(p^\ell+1)(k-1) & =\frac{e_0p^a(p^{2a}-1)}{2(p^{\ell}-1)}-p^\ell (p^{\ell}+1)\\
                      & =\frac{e_0(p^a-1)}{2(p^{\ell}-1)}(p^{2a}+p^a+2)-\frac{e_0(p^a-1)}{p^{\ell}-1}-p^\ell(p^\ell+1),
\end{align*}
 which yields $p^{2a}+p^a+2\mid \frac{e_0(p^a-1)}{p^{\ell}-1}+p^\ell(p^\ell+1)$.
By $p^{2\ell}\leq p^a$, we have $p^{2a}+p^a+2<
 4a(p^a-1)+2p^a$. It follows that $p^{2a}<(4a+1)(p^a-1)-1<(4a+1)p^a$, and then $p^a<4a+1$. This is
 impossible.

\subsubsection*{{\bf Case (5)}. $X\cap G_x= \PSL(2, q_0)$, {\it for} $q=q_0^r$ {\it where} $r$ {\it is an odd prime}.}

Here $v=\frac{q_0^{r-1}(q_0^{2r}-1)}{q_0^2-1}$,
$|X_x|=\frac{1}{2}q_0(q_0^2-1)$, $|\Out(X)|=2f$,
$|G|=\frac{1}{2}eq(q^2-1)$, and $|G_x|=\frac{1}{2}eq_0(q_0^2-1)$,
where $ e\mid  2f$. Let $q_0=p^a$. Then $f=ra$.

From $|G_x|^3>|G|$, that is, $(\frac{1}{2}eq_0(q_0^2-1))^3>\frac{1}{2}eq(q^2-1)=\frac{1}{2}eq_0^r(q_0^{2r}-1)$, we obtain
$$4f^2\geq e^2>4q_0^{r-3}\frac{q_0^{2r}-1}{q_0^6-3q_0^4+3q_0^2-1}.$$
For an odd prime $r$, if $r\geq 5$, then
$$f^2>q_0^{r-3}\frac{q_0^{2r}-1}{q_0^6-3q_0^4+3q_0^2-1}\geq q_0^{r-3} \frac{q_0^{10}-1}{q_0^6-3q_0^4+3q_0^2-1}
>q_0^r=q=p^f,$$
where the third inequality holds because
$q_0^{10}-1>q_0^3(q_0^6-3q_0^4+3q_0^2-1)=q_0^9-3q_0^5(q_0^2-1)-q_0^3$.
But it is easy to see that $\frac{p^f}{f^2}>1$ when  $p\geq 3$ and $f\geq r \geq 5$, a contradiction.
Hence $r=3$, and so $v=q_0^2(q_0^4+q_0^2+1)$ and $f=3a$.

The subdegrees of $\PSL(2,q_0^3)$ on the cosets of $\PSL(2,q_0)$
are (\cite{Subdegree}):
      $$ 1,  \ \big(\frac{q_0^2-1}{2}\big)^{2({q_0}+1)},\ ({q_0}({q_0}-1))^{\frac{{q_0}({q_0}-1)}{2}},
        \ ({q_0}({q_0}+1))^{  \frac{{q_0}({q_0}+1)}{2} },\ \big(\frac{{q_0}(q_0^2-1)}{2}\big)^{2(q_0^3+{q_0}-1)}.$$
 By Lemma
\ref{divide1}, we have
 $$ k \mid   \lambda\gcd\big( ({q_0}+1)^2({q_0}-1), \frac{q_0^2({q_0}-1)^2}{2}, \frac{q_0^2({q_0}+1)^2}{2},
  {q_0}(q_0^2-1)(q_0^3+{q_0}-1) \big).$$
Now, since $ \gcd(q_0+1, q_0-1)= 2$ and $\gcd(q_0,
q_0+1)=\gcd(q_0, q_0-1)=1$, we have
$\gcd\big(\frac{q_0^2({q_0}-1)^2}{2},
\frac{q_0^2({q_0}+1)^2}{2}\big)$
$=q_0^2\gcd\big(\frac{({q_0}-1)^2}{2},
\frac{4q_0}{2}\big)=2q_0^2$,  $\gcd\big( ({q_0}+1)^2({q_0}-1),
2q_0^2 \big)=2$, and so $ k \mid  2\lambda$. Thus $k=2\lam$
follows from $ k > \lambda$. The equation $k(k-1)=\lambda(v-1)$
forces $v=4\lambda-1$. Therefore
$\lambda=\frac{v+1}{4}=\frac{q_0^6+q_0^4+q_0^2+1}{4}$ and
$k=2\lambda=\frac{q_0^6+q_0^4+q_0^2+1}{2}$. Then by Lemma
\ref{divide} ($iii$), $k \mid  |G_x|=\frac{1}{2}eq_0(q_0^2-1)$.
This together with $e\mid2f=6a$ and $q_0=p^a$, implies
$\frac{p^{6a}+p^{4a}+p^{2a}+1}{2}\leq 3ap^a(p^{2a}-1)$ and so that
$ p^{6a}< 6a\cdot p^a\cdot p^{2a}$, i.e., $p^{3a}<6a$, which is
impossible.

\subsubsection*{{\bf Case (6)}. $X\cap G_x= A_5$, {\it where}  $q=p\equiv \pm 1({\rm mod} \ {5})$ {\it  or}   $q=p^2\equiv -1({\rm mod} \ {5})$.}

Here $v=\frac{q(q^2-1)}{120}$,  $|X_x|=|X\cap G_x|=60$,
$|\Out(X)|=2f$, $|G|=\frac{1}{2}eq(q^2-1)$ and $|G_x|=60e$, where
$ e\mid  2f$ and $f=1$ or 2.

From the inequality $|G_x|^3>|G|$ we have
$(60e)^3>\frac{1}{2}eq(q^2-1)$. This together with $e\mid  2f$,
implies $2\cdot 60^3\cdot (2f)^2\geq p^f(p^{2f}-1)$, i.e.,
$$120^3f^2\geq p^f(p^{2f}-1).$$

If $ f=1$ then $q=p\equiv \pm 1\pmod {5}$ and $120^3\geq
p(p^{2}-1)$, which force $q= 11$, 19, 29, 31, 41, 59, 61, 71, 79,
89, 101 or 109. Now we compute the values of $v$ by
$v=\frac{q(q^2-1)}{120}$, and from $k\mid |G_x|$, $e=1$ or 2 we
get $k\mid u=120$. We then check  all possibilities for $v$ by
using Algorithm \ref{algorithm}, and obtain three potential
parameters: $(11, 5, 2)$, $(11, 6, 3)$ and $(57, 8, 1)$.
 If $(v,k,\la)=(57, 8, 1)$, then $X=\PSL(2,19)$.
The subdegrees of $X$ on the cosets of $A_5$ are 1, 6, 20 and 30. 
By Lemma \ref{orbits}, the subdegrees of $G$ are also 1, 6, 20 and
30, contradicting Lemma \ref{divide}($ii$). If $(v,k,\la)=(11, 5,
2)$, then $X=\PSL(2,11)$,  and so $G=\PSL(2,11)$ or $\PGL(2,11)$.
The {\sc GAP}-command  {\tt Transitivity($G,\Omega$)} returns the
degree $t$ of transitivity of the action implied by the arguments;
that is, the largest integer $t$ such that the action is
$t$-transitive. Thus we know that $G$ acts as 2-transitive
permutation group on the set $P$ of 11 points by  {\sc GAP}. Since
$\gcd(k, \lam)=1$, then by Lemma \ref{flag} we see that $\D$ is
flag-transitive, as required. In fact, this design has been found
in \cite{Reg ClassicalBi}. If $(v,k,\la)=(11, 6, 3)$, then Lemma
\ref{ComplementFlagTr} shows that $\D$ is also flag-transitive, as
described in \cite{Zhou ClassicalTri}.

If $f=2$ then $q=p^2\equiv -1\pmod {5}$ and $120^3\cdot 4^2\geq
p^2(p^{4}-1)$. Hence, the possible pairs $(p,v)$ are $(3, 6)$,
$(7, 980)$ and $(13, 40222)$. Since $k\mid  60e$ and $e\mid 2f=4$,
we have  $k\mid u=240$. Running Algorithm \ref{algorithm} with
$u=240$ and $v=6$, 980 or 40222, returns an empty  list {\sc
Designs} for every case, a contradiction.

\subsubsection*{{\bf Case (7)}. $X\cap G_x= A_4$, $q=p\equiv \pm 3\pmod {8}$ {\it and } $ q\not\hspace{-0.2mm}\equiv \pm 1\pmod {10}.$}

Here $v=\frac{q(q^2-1)}{24}$, $|X_x|=|X\cap G_x|=12$,
$|\Out(X)|=2$, $|G|=\frac{1}{2}eq(q^2-1)$ and $|G_x|=12e$, where $\
e=1$ or 2.

The inequality $|G_x|^3>|G|$ gives $(12e)^3>\frac{1}{2}eq(q^2-1)$.
Since $q\geq 5$, $q=p\equiv \pm 3\pmod {8}$  and  $
q\not\hspace{-0.2mm}\equiv \pm 1\pmod {10}$, we get $q=5$ or 13.
Thus $v=5$ or 91, respectively. It is not hard to see that there
is no symmetric $(v,k,\la)$ design with $v=5$. If $v=91$ then all
possible parameters of $(k,\lam)$ are
$$ (10, 1), (36, 14), (45, 22), (46, 23),  (55, 33)\  \hbox{and}\ (81, 72).$$
However, by $k\mid 12e$ and $e=1$ or 2, we have $k\mid 24$, the
desired contradiction.

\subsubsection*{{\bf Case (8)}. $X\cap G_x= S_4$, $q=p\equiv \pm 1\pmod {8}.$}

Now $v=\frac{q(q^2-1)}{48}$, $|X_x|=|X\cap G_x|=24$, $|\Out(X)|=2$,
$|G|=\frac{1}{2}eq(q^2-1)$, $|G_x|=24e$, where $e=1$ or 2.

Since $q=p$, $e\leq 2$ and $|G_x|^3>|G|$, that is, $(24e)^3>\frac{1}{2}eq(q^2-1)$,
we  get
$$q(q^2-1)<2\cdot 24^3\cdot e^2\leq 48^3.$$
Since $q\equiv \pm 1\pmod {8}$, we obtain that the possible pairs $(q,v)$ are
(7, 7), (17, 102), (23, 253), (31, 620), (41, 1435) and (47, 2162).
 Since $k\mid  |G_x|=24e$ and $e=1$ or 2, we get $k\mid u=48$.
Thus Algorithm \ref{algorithm} gives only two parameters: $(7, 3,
1)$ and $(7, 4, 2)$. If $(v,k,\la)=(7, 3, 1)$, then $X=\PSL(2,7)$,
and so $G=\PSL(2,7)$ or $\PGL(2,7)$. Hence  $G$ acts as a
2-transitive permutation group on the set $P$ of 7 points by  {\sc
GAP}.
 Lemma \ref{flag} shows that $\D$ is flag-transitive because $\gcd(k, \lambda)=1$.
If $(v,k,\la)=(7, 4, 2)$, then $\D$ is also flag-transitive by
Lemma \ref{ComplementFlagTr}. This design has been discussed in
\cite{Reg ClassicalBi}.

\subsection{Characteristic two}

In this subsection, we suppose that $G$ is of characteristic 2 and $X\cap G_x$ is maximal in $X$.
The structure of $X\cap G_x$  is given in Table  \ref{tab:2}.

\subsubsection*{{\bf Case (1)}. $X\cap G_x= E_q: (q-1).$}

Here $v=q+1$ and $k(k-1)=\lam (v-1)=\lam q=2^f\lam $. If $2\mid
k$, then $\gcd(2,k-1)=1$, and so $2^f\mid  k$, which implies $v-1
\mid  k$. This is impossible. Thus $2\nmid k$ and so $ 2^f\mid
k-1$. It follows that $v-1\mid  k-1$,  a contradiction.

\subsubsection*{{\bf Case (2)}. $X\cap G_x= D_{2(q-1)}.$ }

Now $v=\frac{1}{2}q(q+1)$,  $|\Out(X)|=f$, $|G|=eq(q^2-1)$ and
$|G_x|=2e(q-1)$, where $e\mid  f$.

From $q=2^f\geq 4$  we know that  $v=\frac{1}{2}q(q+1)$ is even.
So $\la$ is also even since $k(k-1)=\lam (v-1)$. Lemma
\ref{divide} ($iii$) shows $k\mid  2e(q-1)$. Then there exists a
positive integer $m$ such that $k=\frac{2e(q-1)}{m}$. Again by
$k(k-1)=\lam (v-1)$, we have
$\frac{2e(q-1)}{m}(\frac{2e(q-1)}{m}-1)=\lam
(\frac{1}{2}q(q+1)-1)$, and so
$(8e^2-m^2\lam)q=2m^2\lam+8e^2+4em$, which forces $8e^2-m^2\lam>0$
and so $m<2e$. The fact that $\la$ is even implies that
$8e^2-m^2\lam\geq 2$. So we have
 $$ 2^f=q=\frac{24e^2+4em}{8e^2-m^2\lam}-2\leq \frac{24e^2+4e\cdot 2e}{2}\leq 16f^2.$$
 Hence $2\leq f\leq 10$.
Since $k\mid 2e(q-1)$ and $e\mid  f$, we get $k\mid u=2f(q-1)$.
The pairs $(v,u)$, for $2\leq f\leq 10$, are (10, 12), (36, 42),
(136, 120), (528, 310),
  (2080, 756), (8256, 1778), (32896, 4080), (131328, 9198) and
  (524800, 20460).
Then Algorithm \ref{algorithm} gives only one possible set of
parameters $(36, 21, 12)$. Suppose $(v,k,\la)=(36, 21, 12)$. Then
$G=\PSL(2,8)$ or $\PGGL(2,8)$. When $G=\PSL(2,8)$, the subdegrees
of $G$  are 1, $7^3$ and 14, and $G$ has only one conjugacy class
of subgroups of index 36. Thus for any $B\in \B$, $G_x$ is
conjugate to $G_B$. Without loss of generality, let $G_x=G_{B_0}$
for some block $B_0$. The flag-transitivity of $G$ forces
$G_{B_0}$ to act transitively on the points of $B_0$. Hence the
points of $B_0$ form an orbit of $G_x$, which implies that a
subdegree of $G$ is $k=21$, a contradiction. Now assume
$G=\PGGL(2,8)$. Then  the subdegrees of $G$ are 1, 14 and 21, and
$G$ has only one conjugacy class of subgroups of index 36. So let
$G_x=G_{B_0}$ for some block $B_0$ as above. Then $B_0$ is an
orbit of size 21 of $G_x$. By using {\sc Magma}, we obtain that
$|\B|=|B^G|=36$, but $|B_i\cap B_j|=10$ or 15 for any two distinct
blocks $B_i$ and $B_j$. This is a contradiction since in our
situation any two distinct blocks should have $\lambda=12$ common
points.

\subsubsection*{{\bf Case (3)}.  $X\cap G_x= D_{2(q+1)}.$}

Here $v=\frac{1}{2}q(q-1)$, $|\Out(X)|=f$,
$|G|=\frac{1}{2}eq(q^2-1)$ and $|G_x|=2e(q+1)$, where $e\mid  f$.

Since $k\mid  |G_x|$, there exists a positive integer $m$ such
that $k=\frac{2e(q+1)}{m}$. Thus Lemma \ref{divide} ($i$) yields
$\frac{2e(q+1)}{m}\big(\frac{2e(q+1)}{m}-1\big)=\lam
(\frac{1}{2}q(q-1)-1)$, and so
$(m^2\lam-8e^2)q=8e^2-4em+2m^2\lam=8(e-\frac{1}{2}m)^2+2(\lam-1)m^2>0$.
We then have
$$ 2^f=q=\frac{8e^2-4em+2m^2\lam}{m^2\lam-8e^2}=\frac{24e^2-4em}{m^2\lam-2e^2}+2,$$
which implies $2^f<24e^2+2\leq 24f^2+2.$ Hence $2\leq f\leq 11$.
Since $k\mid 2e(q+1)$ and $e\mid  f$, we have  $k\mid u=2f(q+1)$.
For $2\leq f\leq 11$, the pairs $(v,u)$  are as follows:
$$
\begin{array}{lllllll}
  &(6, 20), &(28, 54), &(120, 136), &(496, 330), &(2016, 780), \\
  &(8128, 1806), &(32640, 4112), &(130816, 9234), &(523776, 20500), &(2096128,
   45078).&
\end{array}
$$

Applying Algorithm \ref{algorithm} to these pairs $(v,u)$, we
obtain  $(v,k,\la)=(496, 55, 6)$ or $(2016$, $156$, $12)$. If
$(v,k,\la)=(496, 55, 6)$, then $G=\PSL(2, 2^5)$ or $\PGGL(2,
2^5)$. Let $G=\PSL(2, 2^5)$(or $\PGGL(2, 2^5)$), then the
subdegrees of $G$ are 1 and $33^{15}$ (or 1 and $165^3$),
 and $G$ has only one conjugacy class of subgroups of index 496.
Thus there exists a block-stabilizer $G_{B_0}$ such that
$G_x=G_{B_0}$, which implies that $B_0$ should be an orbit of
$G_x$. But this is impossible because $|B_0|=55$. Now suppose
$(v,k,\la)=(2016, 156, 12)$. Then $G=\PSL(2, 2^6)$, $\PSL(2,
2^6):i$ $(i=2,3)$ or $\PSSL(2, 2^6)$. By the fact that $G$ has
only one conjugacy class of subgroups of index 2016,
 similar to the analysis above,  there exists a block $B_0$ such that $B_0$ is an orbit of $G_x$.
 Thus $G_x$ should have an orbit of size 156.
The subdegrees of $G$, however, are as follows:
\begin{enumerate}
  \item[(i)] \, 1, and $65^{31}$ when  $G=\PSL(2, 2^6)$;

  \item[(ii)] \, 1, $65^7$ and $130^{12}$ when  $G=\PSL(2,
  2^6):2$;

  \item[(iii)] \, 1, $65$ and $195^{10}$ when $G=\PSL(2, 2^6):3$;
  \item[(iv)] \, 1, 65, $195^2$ and $390^4$ when $G=\PSSL(2,
2^6)$.
\end{enumerate}

\subsubsection*{{\bf Case (4)}.  $X\cap G_x=\PSL(2, q_0)=\PGL(2, q_0)$, {\it where} $q=q_0^r$
 {\it  for some prime} $r$ {\it  and} $q_0\neq 2$.}

 Here $v=\frac{q_0^{r-1}(q_0^{2r}-1)}{q_0^2-1}$, $|X_x|=|X\cap G_x|=q_0(q_0^2-1)$,   $|\Out(X)|=f$,
$|G|=\frac{1}{2}eq(q^2-1)$ and $|G_x|=eq_0(q_0^2-1)$, where $e\mid
f$. Let $q_0=2^a$, so that $f=ra$.

  From $|G_x|^3>|G|$, $q=q_0^r$ and $e\mid f$, we get
$$f^2\geq e^2>q_0^{r-3}\frac{q_0^{2r}-1}{q_0^6-3q_0^4+3q_0^2-1}.$$
If $r\geq 5$, then
$$f^2>q_0^{r-3}\frac{q_0^{2r}-1}{q_0^6-3q_0^4+3q_0^2-1}\geq q_0^{r-3} \frac{q_0^{10}-1}{q_0^6-3q_0^4+3q_0^2-1}
>q_0^r=q=2^f,$$
where the third inequality holds because
$q_0^{10}-1>q_0^3(q_0^6-3q_0^4+3q_0^2-1)=q_0^9-3q_0^5(q_0^2-1)-q_0^3$.
But for $f\geq r \geq 5$ the inequality $f^2>2^f$ is not satisfied.
 Hence $r=2$ or 3.

Suppose first that $r=3$, so that $q=q_0^3=2^{3a}$,
$v=q_0^2(q_0^4+q_0^2+1)$ and $f=3a$. The subdegrees of
$\PSL(2,q_0^3)$ on the cosets of $\PSL(2,q_0)$ are as follows
(\cite{Subdegree}):
      $$ 1,  \ (q_0^2-1)^{{q_0}+1},\ ({q_0}({q_0}-1))^{\frac{{q_0}({q_0}-1)}{2}},
        \ ({q_0}({q_0}+1))^{  \frac{{q_0}({q_0}+1)}{2} },\ ({q_0}(q_0^2-1))^{q_0^3+{q_0}-1}.$$
By Lemma \ref{divide1}, we have
 $$ k \mid   \lambda\gcd\big( ({q_0}+1)^2({q_0}-1), \frac{q_0^2({q_0}-1)^2}{2}, \frac{q_0^2({q_0}+1)^2}{2},
  {q_0}(q_0^2-1)(q_0^3+{q_0}-1) \big).$$
 So $ k \mid  2\lambda$.
 This forces $k=2\lam$ since $k > \lambda$.
Thus  $v=4\lambda-1$ by equation $k(k-1)=\la (v-1)$. Then
$\lambda=\frac{v+1}{4}=\frac{q_0^6+q_0^4+q_0^2+1}{4}$ and
$k=2\lambda=\frac{q_0^6+q_0^4+q_0^2+1}{2}$. By $k \mid
|G_x|=\frac{1}{2}eq_0(q_0^2-1)$ and $e\mid f=3a$,
 we get $\frac{2^{6a}+2^{4a}+2^{2a}+1}{2}\leq \frac{3a}{2}\cdot
 2^a(2^{2a}-1)$,
and so $ 2^{6a}\leq 3a\cdot 2^a\cdot 2^{2a}$, i.e., $2^{3a}\leq
3a$, which is impossible.

Now suppose $r=2$. Then $q=q_0^2=2^{2a}$, $v=q_0(q_0^2+1)$ and
$f=2a$. The subdegrees of $\PSL(2,q_0^2)$ on the cosets of
$\PGL(2,q_0)$ are as follows (\cite{Subdegree}):
      $$ 1,\ q_0^2-1,\ (q_0(q_0-1))^{\frac{q_0-2}{2}},\ (q_0(q_0+1))^{\frac{q_0}{2}}.$$
By Lemma \ref{divide1}, we have
$$k \mid  \lambda\gcd\big(q_0^2-1, \frac{q_0(q_0-1)(q_0-2)}{2}, \frac{q_0^2(q_0+1)}{2}\big).$$
Here
   $\gcd\big(\frac{q_0(q_0-1)(q_0-2)}{2}, \frac{q_0^2(q_0+1)}{2}\big)
        =2^{a}\gcd((2^a-1)(2^{a-1}-1), 2^a+1)$
which divides  $3\cdot 2^{a}=3q_0$. Thus $ k\mid
\lambda\gcd\big(q_0^2-1, 3q_0\big)$, and so $k\mid 3\lambda$. Now,
$k > \lambda$ implies that $k = 3\lambda$ or $\frac{3\lambda}{2}$.

If $k = 3\lambda$, then $v=9\lambda-2$ by $k(k-1)=\la (v-1)$. So
$\lambda=\frac{v+2}{9}=\frac{q_0^3+q_0+2}{9}$ and  $
k=3\lambda=\frac{q_0^3+q_0+2}{3}$. From $k \mid
|G_x|=eq_0(q_0^2-1)$ and $e\mid f=2a$, we have $k \mid
2aq_0(q_0^2-1)$. By the facts that $\gcd(q_0^3+q_0+2, q_0)=2$ and
$\gcd(q_0^3+q_0+2, q_0-1)=\gcd(4, q_0-1)=1$, we get
$\frac{q_0^2-q_0+2}{3}\mid  4a$, and so $\frac{2^{2a}-2^{a}+2}{3}
\leq  4a$, which implies that $a=1$ or 2. Since $q_0\neq 2$,
$a\neq 1$. Hence $a=2$ and $q_0=4$, but then $k=\frac{70}{3}$ is
not an integer.

If $k=\frac{3\lambda}{2}$, then $v=\frac{9\lambda-2}{4}$. Thus
$\lambda=\frac{4v+2}{9}=\frac{4q_0^3+4q_0+2}{9}$ and
$k=\frac{2q_0^3+2q_0+1}{3}$. Since $k \mid  |G_x|$ and $e\mid
f=2a$, we have $\frac{2q_0^3+2q_0+1}{3} \mid  2aq_0(q_0^2-1)$.
Note that $\gcd(2q_0^3+2q_0+1, q_0)=1$, $\gcd(2q_0^3+2q_0+1,
q_0+1)=1$ or 3, and $\gcd(2q_0^3+2q_0+1, q_0-1)=1$ or 5. Then
$2q_0^3+2q_0+1\mid 90a$, and hence $2^{3a+1}+ 2^{a+1}+1\leq 90a$.
It follows that $a=1$ or 2. If $a=1$ then $q_0=2$, a
contradiction. If $a=2$ then $q_0=4$ which implies
$k=\frac{137}{3}$ is not an integer.

This completes the proof of Theorem 1.1.  $\hfill\square$


\end{document}